\documentclass[11pt]{article}
\usepackage{a4wide}
\usepackage{amsmath}

\usepackage{amsfonts}
\usepackage{amssymb}
\usepackage{euscript}
\usepackage[latin1]{inputenc}
\usepackage{graphicx}

\newenvironment{proof}{\noindent {\it Proof.~~}\ }{\  \rule{1mm}{2mm}\medskip}

\newenvironment{proof*}{\noindent {\it Proof.~~}\ }{}

\newtheorem{theorem}{Theorem}
\newtheorem{lemma}[theorem]{Lemma}
\newtheorem{corollary}[theorem]{Corollary}

\newtheorem{theirtheorem}{Theorem}

\newtheorem{theirlemma}[theirtheorem]{Lemma}

\usepackage{latexsym,amssymb}
\usepackage{amsmath}

\newcommand{\subgp}[1]{\langle{#1}\rangle}

\def\qed{\hfill$\Box$}

\begin{document}
\title{On vosperian and superconnected vertex-transitive digraphs}

\author{ {Y. O. Hamidoune}\thanks{
Universit\'e Pierre et Marie Curie, E. Combinatoire, Case 189, 4
Place Jussieu, 75005 Paris, France. \texttt{yha@ccr.jussieu.fr}}
\and {A. Llad\'o}\thanks{
 Universitat
Polit\`ecnica de Catalunya, Dept. Matem\`atica Apl. IV;  Jordi
Girona, 1, E-08034 Barcelona, Spain. \texttt{allado@ma4.upc.edu}}
\and {S. C. L\'opez}\thanks{
 Universitat
Polit\`ecnica de Catalunya, Dept. Matem\`atica Apl. IV;  Jordi
Girona, 1, E-08034 Barcelona, Spain. \texttt{susana@ma4.upc.edu}} }

\date{}

\maketitle

\begin{abstract} We investigate the structure of a digraph  having a transitive automorphism group where every cutset of minimal cardinality consists of
 all successors or all predecessors of some vertex. We improve most of the existing results in this area.
\end{abstract}

\section{Introduction}

By a digraph, we shall mean a finite directed graph  having no loops. Let   $\Gamma =(V,E)$  be a digraph.
The set successors (resp. antecedents) of a vertex  $x\in V$ will be denoted by $\Gamma (x)$
( $\Gamma^- (x)$ ). The degree of $x$ is defined as $|\Gamma (x)|.$ If all the vertices have the same degree, the graph will be called regular
and this common value will be called the degree of $\Gamma.$
A subset $C$ of  $V$ is said to be {\em  strongly connected} if, for any two vertices $x,y$ of $C$ there is a directed
path of $\Gamma$ from $x$ to $y$ with all vertices contained in $C$.
The digraph $\Gamma$ is said to be {\em  strongly connected} if  $V$ is  strongly connected.
A subset $T$ of $V$ will be called a  {\em cutset} of $\Gamma$ if $V\setminus T$ is not  strongly connected. The minimal cardinality of a cutset   is called the {\em
 connectivity} of ${\Gamma}$. %The {\em degree} of a vertex $v\in V$ is the cardinality of

Regular digraphs may have small connectivity. But vertex transitive digraphs have a large one. It was proved independently by
Mader \cite{mader} and Watkins \cite{watkins} that the connectivity of a  strongly connected vertex-transitive symmetric digraph with degree $d$
is greater than $2d/3.$ One of the authors proved in \cite{Ham77} that a  strongly connected vertex-transitive  digraph with degree $d$ has connectivity greater $ d/2.$
It was proved by Watkins  in \cite{watkins} that the connectivity of a  strongly connected edge-transitive graph is its minimal degree.
One of the authors proved in \cite{Ham77} that the connectivity of a  strongly connected arc-transitive digraphs is its  degree.

 A digraph $\Gamma$ with degree $d$ is said to be {\em
 superconnected} if every cutset with cardinality less than $d+1$   consists of
 all successors or all predecessors of some vertex. This notion,
 introduced in the undirected case, by Boesch and Tindell
 \cite{BoeTin84}, was  also investigated among others by Fiol
 \cite{Fio94} and Balbuena and Carmona \cite{BalCar01}. For an introduction to superconnectivity and related topics, we recommend  the
 chapter by F\`abrega and Fiol in \cite{FabFio04}.

 If any   cutset $T$  with  $|T|\le d$ and $|T|\neq  |V|-3,$ creates exactly two  strongly connected components one of them
 consists of a single vertex, the digraph $\Gamma$ will be called
 {\em vosperian}.
 Clearly  {
 vosperian} digraphs are superconnected.
  This notion was introduced by the authors of
 \cite{HamLlaSer91} as a step in their characterization of
 superconnected Abelian Cayley digraphs.  As observed by Meng in \cite{Men03},  vosperianity is basically equivalent to
 the notion of hyper-connectdness, introduced independently in the undirected case, by Boesch \cite{Boe86}.

A pair of vertices  $\{x,y\}$ will be called a {\em twin pair} (resp. {\em anti-twin pair}) if $x$ and $y$ have the same successors (resp. predecessors ). Following Meng \cite{Men03}, we shall say that a digraph having a twin pair or an  anti-twin pair is {\em reducible}. The presence
of twins or anti-twins is clearly an obstruction for vosperianity.

A recursive characterization of vosperian and superconnected Abelian
Cayley digraphs was obtained by the authors of \cite{HamLlaSer91}.
One of the authors \cite{Ham97} obtained a non-recursive characterization using
the $2$-atoms.

The description of vosperian and superconnected vertex-transitive
digraphs has  received considerable attention in recent years. Most
of the characterizations use a result due to one of the authors in
\cite{Ham97}. The case of edge-transitive undirected graphs was
considered by Meng \cite{Men03}. Liang, Meng and Zhang
\cite{LiaMenZha07} described the case of bipartite
 undirected graphs. Arc transitive digraphs are investigated by Zhang and
Meng in \cite{MenZha09}. In the last paper it is shown that a
strongly connected irreducible anti-symmetric  arc-transitive digraph is
either vosperian or bi-star superconnected. The last notion is
avoidable in our approach, but the reader may refer to
\cite{MenZha09} for its definition. It is also shown that strongly connected
irreducible anti-symmetric arc-transitive digraphs are
superconnected.

More recently, the authors of \cite{HamLlaLop} have given necessary
and sufficient conditions for an undirected vertex-transitive graph
to be vosperian. They have also proved that an irreducible
superconnected vertex transitive undirected graph is vosperian
\cite{HamLlaLop}.

We investigate in the present paper, the vosperianity and
superconnectdness of vertex-transitive digraphs. We give a complete
characterization of vosperian arc-transitive digraphs (Theorem 4)
which completes the result of \cite{MenZha09}. It states that an
arc-transitive   strongly connected digraph is vosperian if and only if it is
irreducible. In particular, this is the case if the degree is
coprime with the order of the digraph. We give also a complete
characterization of vosperian  Cayley digraphs (Theorem 6) and  a
complete characterization of irreducible superconnected Cayley
digraphs (Theorem 7). These two last characterizations extend the
corresponding ones in the directed case in \cite{HamLlaSer91,Ham97},
and the ones for the undirected case in \cite{HamLlaLop}.

The paper is organized as follows. After giving some terminology and
basic results in Section 2 we describe our main tools in Section 3.
Section 4 is devoted to the study of the arc--transitive case. We
then focus on Cayley digraphs and give the characterization of
vosperianity in Section 5 and of superconnectedness for Cayley
digraphs without twin pairs (i.e. defined by an aperiodic subset) in
Section 6.

%The  organization of the paper is as follows:

% Sections 2,  and 3 contain  our terminology and known tools.
% In Section 4,  we show that an arc-transitive digraph without twin or anti-twins is vosperian.

% In Section 5,  we give a complete  characterization of vosperian  Cayley
%digraphs. In Section 6, we give a characterization of
%superconnected Cayley digraphs without twins or anti-twins (i.e. defined by an aperiodic subset).

\section{Terminology and preliminaries}

Let $\Gamma =(V,E)$ be a
digraph.
The elements of $V$ will be called {\em vertices}. The
elements of $E$ will be called {\em arcs}. Undirected graphs  can be safely identified with symmetric
digraphs. The {\em reverse} of
$\Gamma$ is the digraph  $\Gamma^{-}$ obtained by reversing the
orientation of the arcs of $E.$
% The set of successors of a vertex $x$ will be denoted by $\Gamma (x)$. The set of predecessors of a vertex $x$  will be denoted by $\Gamma ^-(x)$.
More formally, $\Gamma^{-}=(V,E^-),$ where $E^-=\{(x,y) : (y,x)\in E\}.$

%\begin{itemize}
 % \item $\Gamma^{-}=(V,E^-),$ where $E^-=\{(x,y) : (y,x)\in E\},$
  %\item $\Gamma (x)=\{y : (x,y)\in E\},$
  %\item $\Gamma ^- (x)=\{y : (y,x)\in E\}.$
%\end{itemize}

For a regular digraph $\Gamma,$ we denote the degree of any vertex by
$d(\Gamma).$
For a set $X\subset V$  we write $\Gamma (X)=\bigcup_{x\in X}\Gamma
(x)$  and
$\Gamma [X]$ will denote the subdigraph induced by $X$.
A {\em source} of $\Gamma$ is a  subset  $X\subset V$
such that $\Gamma^- (X)\subset X.$ A {\em sink} of $\Gamma$ is a
subset  $X\subset V$ such that $\Gamma (X)\subset X.$  It is
well known that a finite digraph is  strongly connected if and only
if it has no a proper subset which is a source (sink).
As a warning we mention that our $\Gamma (X)$ is written $\Gamma^+ (X)$
in some Graph Theory text books. Also our degree is sometimes called outdegree.

An {\em automorphism } of $\Gamma $  is a bijection $f
:V\rightarrow V$ such that $f (\Gamma (x))=\Gamma(f(x)),$  for every
$x\in V.$ A digraph is said to be {\em vertex-transitive} if, for
any pair $x,y$ of vertices, there is an automorphism that maps $x$
into $y$. A digraph is said to be {\em arc-transitive} if, for
any pair $(x,y),(x',y')$ of arcs, there is an automorphism that maps $x$
into $x'$ and maps $y$
into $y'$. It is an easy exercise to show that a  strongly connected arc-transitive digraph is also vertex-transitive,
and hence regular.

Let $G$ be a group and let $S\subset G\setminus\{1\}.$ The
digraph $\mbox{Cay} (G,S)=(G,E),$ where $E=\{(x,y): x^{-1}y\in S\}$
is called a {\it Cayley} digraph.
Recall that  $\Gamma$ is  strongly connected  if and only if  $S$ is a generating subset.
Note that the  left translation $
\gamma_a : x \mapsto ax$ is an automorphism of $\mbox{Cay} (G,S).$
A set $S$ of a group $G$ is said to be {\em left-aperiodic} if for some $x\in
G\setminus \{1\}$, $xS=S$. If both  $S$ and $S^{-1}$ are not left-periodic, we shall say that
$S$ is {\em aperiodic}.
The next lemma is just an exercise:

\begin{lemma} Let  $S$ be a subset of a group $G.$  Then
\begin{itemize}
  \item[(i)]  $\mbox{Cay} (G,S)=(G,E)$ has a pair of twins if and only if $S$ is left-periodic,
  \item[(ii)]  $\mbox{Cay} (G,S)=(G,E)$ is irreducible  if and only if $S$ is aperiodic.
\end{itemize}
\label{twincayley}
 \end{lemma}

%\subsection{Isoperimetric operators}

%The minimal degree is $\delta_\Gamma =\min\{|\Gamma (x)|: x\in V\}$.
Let $\Gamma =(V,E)$ be digraph.
 Given a subset $X\subset V,$ the {\em boundary}  of $X$ is defined as
  $$\partial_{\Gamma} (X)=\Gamma (X)\setminus X.$$

The {\em exterior} of $X$ is defined as
$
\nabla_{\Gamma} (X)=V\setminus (X\cup \Gamma (X)).
$
We write $\partial^-_{\Gamma}(X)=\partial_{\Gamma^-}(X)$ and  $
\nabla^-_{\Gamma}(X)=\nabla_{\Gamma^-}(X).
$
If the context is clear the reference to $\Gamma$ will be omitted.
Every set $X$ induces the partition (with possibly empty parts) $\{ X,\partial (X), \nabla (X)\}$ of
the vertex set with no arc from $X$ to $\nabla (X)$.

In particular,
\begin{equation}\label{boarduality}
\partial^- (\nabla (X))\subset \partial (X).
\end{equation}

The digraph $\Gamma$ is said to be $k$-{\em separable} if there is a subset $X\subset V$ such that
$\min\{|X|,|\nabla (X)|\}\ge k$.
We shall say that a subset $X$ of $V$ induces a $k$-separation  on $\Gamma$ if $k\le \min\{|X|,|\nabla (X)|\}< \infty$. In particular, $X$ induces a  $k$-separation on $\Gamma$ if and only if $\nabla (X)$ induces a  $k$-separation on $\Gamma^{-}$. Hence, $\Gamma$ is $k$-separable if and only if $\Gamma^{-}$ is  $k$-separable.

%\subsection{Cutsets}
As an exercise, the reader may check that a subset $T\subset V$ is a { cutset} if and only if there is a nonempty subset $X\subset V,$ with  $\nabla (X)\neq \emptyset$
and $\partial (X)\subset T.$ One may see from (\ref{boarduality}) that  $T$ is a cutset of $\Gamma$ if and only if
$T$ is a cutset of $\Gamma ^-.$ A cutset of minimal cardinality, where $\min (\emptyset)=|V|-1,$ will be called a
{\em minimum} cutset.

%The {\em connectivity} $\kappa (\Gamma )$  of the digraph $\Gamma$ is traditionally  the cardinality of a minimum cutset.

%Let $\Gamma$ be a finite $d$-regular digraph with $\kappa =d.$ We shall say that  $\Gamma$ is {\em superconnected},
%if for every minimum cutset $T,$ there exists $v\in V\setminus T$ such that
%$\Gamma (v)=T$ or $\Gamma ^- (v)=T.$

For a $k$-separable digraph  $\Gamma,$ the
$kth$-{\em isoperimetric connectivity}  of $\Gamma$ is defined as
$$
\kappa_k (\Gamma)=\min \{|\partial (X)|: X\subset V \ \text{and} \ \min\{|X|,|\nabla (X)|\}\ge
k\}.
$$

A subset $F\subset V$ is called a $k$-{\em fragment} if $\min\{|F|,|\nabla (F)|\}\ge
k$ and $|\partial (F)|=\kappa_k (\Gamma )$. A $k$-{\em atom} is a $k$-fragment
of minimum cardinality.
A $k$-fragment (a $k$-atom) of $\Gamma^{-}$ is called a {\em  negative} $k$-{\em fragment} ({\em a negative $k$-atom}). We also write
$\kappa _{-k}(\Gamma)=\kappa _k(\Gamma^{-})$.

The notion of $kth$-{isoperimetric connectivity} was introduced in \cite{Ham96}.
We recall the following duality lemma:

\begin{theirlemma}[\cite{Ham08}]
Let $\Gamma =(V,E)$ be a  finite $k$-separable digraph. Then $\kappa_k(\Gamma )=\kappa_{-k}(\Gamma )$.\label{equal_kappa}
If $X$ is a $k$-fragment,  then
$$
   \partial^{-}(\nabla  (X))=\partial (X)
\ \text{and} \
   \nabla^-(\nabla (X))=X.$$

In particular, $\nabla(X)$ is a negative $k$-fragment.
\end{theirlemma}

It follows easily that $\kappa(\Gamma)=\kappa_1(\Gamma)$.
We need a special case $k=2$ of the easy next lemma:
\begin{theirlemma}[folklore]\label{katom}
Let $\Gamma =(V,E)$ be a  finite $k$-separable digraph and let $A$ be a $k$-atom with $|A|>k.$
 Then $\Gamma ^-(x)\cap A\neq \emptyset,$ for every $x\in A.$ In particular,  $A$ contains a circuit.
\end{theirlemma}
%\subsection{$2$--atoms}

\section{Some tools}

%\subsection{Minkowski product}
Let $G$ be a group and let $S$ be a subset of $G$. The subgroup
generated by $S$ will be denoted by $\subgp{S}$. Let  $ A,B$ be
subsets of $ G $. The {\em Minkowski product} is defined as
$$AB=\{xy \ : \ x\in A\  \mbox{and}\ y\in
  B\}.$$

We use the following easy lemmas:

\begin{theirlemma}[\cite{manlivre}, Theorem 1]
Let $G$ be a finite group and let
 $A,B$  be  subsets of $G$
 such that $|A|+|B|>|G|$.
 Then $AB=G$.

\label{prehistorical}
 \end{theirlemma}

\begin{theirlemma}[folklore]
Let $G$ be a finite cyclic group generated by $r$ and let
 $B=\{1,r, , r^{|B|-1}\}b.$ For every subset $A\subset G,$
 such that $|AB|=|A|+|B|-1<|G|$, there is an $a\in G$ such that
$aA=\{1,r,\ldots  , r^{|A|-1}\}.$

\label{AP}
 \end{theirlemma}

Let $\Gamma =(V,E)$ be a finite digraph.
 A {\em block of imprimitivity} (or simply a block in what follows) of $\Gamma$ is
 a subset $B\subset V$ such that for every automorphism $f$ of $\Gamma$, either $f(B)=B$ or $f(B)\cap B=\emptyset$.

 Recall the following result:

\begin{theirlemma} [folklore]
Let $G$ be a finite group and let
 $S$  be  a subset of $G.$ Any block $B$ of  $\Gamma =Cay(G,S)$ with $1\in B$
 is a subgroup of $G.$
\label{bloccayley}
 \end{theirlemma}

%\begin{theirtheorem}[\cite{Ham84}]
%Let $\Gamma =(V,E)$ be a  finite  Cayley digraph. Let $H$ be  a $1$-atom containing $1$ and let $K$
%be a negative $1$-atom. If $|H|\le |K|,$ then $H$ is a subgroup.
%\label{atomsg}

%\end{theirtheorem}

As an easy exercise, the reader may show that a digraph $\Gamma$ is vosperian   if and only
if it is either non $2$-separable or if $\kappa _2(\Gamma )\ge d(\Gamma)+1$.
Thus, for finite digraphs Lemma \ref{equal_kappa}  implies that $\Gamma$  is non-vosperian if and only if $\Gamma^-$ is non-vosperian.

We shall use some results from \cite{Ham99}. These results were originally formulated
for directed graphs with loops.
 But they apply to our context, since the deletion (or addition) of loops does not modify the notions of $k$-fragments
and $k$-isoperimetric connectivity. We use the following result:

\begin{theirtheorem} [\cite{Ham99}] Let $\Gamma =(V,E)$ be a  finite non-vosperian vertex-transitive digraph.  Then one of the following holds:
\begin{itemize}
\item[(i)]  There is a $2$-atom of size $2$ or a negative $2$-atom of size $2.$
 \item[(ii)] There is a block which is a $2$-fragment or a negative $2$-fragment.
  \item[(iii)] Every vertex of $V$ is contained in at most two distinct $2$-atoms. Moreover the intersection of two distinct $2$-atoms has cardinality $<2.$
  \item[(iv)] Every vertex of $V$ is contained in at most two distinct negative $2$-atoms.
  Moreover the intersection of two distinct negative $2$-atoms has cardinality $<2.$
\end{itemize}
\label{asterix} \end{theirtheorem}

Note that the notion of superatoms, the smallest atoms with cardinality larger than one, is used in \cite{Ham99} instead of the close
notion of $2$-atoms used here.
%Notice that the paper \cite{Ham99} uses the obsolete notion of superatoms. A simplified proof of the last result using the $2$-atoms terminology is contained in \cite{Ham10}.
%In any case, Theorem \ref{asterix} is the combination of Corollary 2.8, Proposition 4.4, Proposition 4.5 and Proposition 4.8 of \cite{Ham99}.
%We give a small hint below:

%Let $M$ be a $2$-atom
%and  let $N$ be a negative $2$-atom.
%If $\kappa _2(\Gamma)< d(\Gamma)$, then $M$ and $N$ are respectively a $1$-atom and a negative $1$-atom. By  Corollary 2.8, $M$ is a block or $\nabla (M)$ is a block.
%Assume now  $\kappa _2(\Gamma)=d(\Gamma).$  If $|\nabla (M)|<|M|,$ then  $\nabla (M)$ is a block by Proposition 4.4.
%Assume now that $|\nabla (M)|\ge |M|,$ and that $|\nabla (N)|\ge |N|.$ By Proposition 4.5, the intersection of two (rep. negative) distinct $2$-atoms has cardinality $<2.$ Also, we have by Proposition 4.8, that provided (i) and (ii) do not hold then (iii) holds or (iv) holds.

\begin{corollary}[\cite{Ham99}]
Let $S$ be a  generating subset of a finite  group $G$ with $1\notin S$ such that
 $\Gamma =Cay(G,S)$ is non-vosperian.
Then there are  a subgroup $H$ and an  $a\in G$
such that  $H\cup Ha$ is a $2$-fragment or a negative $2$-fragment.
\label{vospercorog}
\end{corollary}

\proof
Let $M$ be a $2$-atom
and  let $N$ be a negative $2$-atom such that $1\in M\cap N.$ We may assume that $\min (|M|,|N|)\ge 3.$ Otherwise the result holds with $H=\{1\}.$ Moreover,  $\Gamma$ has no block which is a $2$-fragment or a negative $2$-fragment.
Otherwise the result holds by Lemma \ref{bloccayley}.

By Theorem \ref{asterix},  $1$ is contained in at most two $2$-atoms of $\Gamma$ or two $2$-atoms of $\Gamma^-.$
Up to replacing  $\Gamma$ by $\Gamma^-,$ we may assume that  $1$ is contained in at most two $2$-atoms of $\Gamma.$
Let
$H=\{x\in G: \hspace{0.2cm} xM=M\}$.
The result holds clearly if $M=H.$ Take an $a\in M\setminus H.$ We shall show that $M=H\cup Ha$.
Observe that $M$ and $a^{-1}M$ are two distinct $2$-atoms containing $1.$ Thus  for every $x\in M,$
we have $x^{-1}M=M$ or $x^{-1}M=a^{-1}M.$ Therefore either $x\in H$ or $xa^{-1}\in H$. In particular, $M\subset H\cup Ha.$ But $HM\subset M,$ by the definition of $H$.\qed

 Let $G$ be an Abelian group and let  $1\in A$ is a  $2$-atom of $Cay(G,S).$ It was proved in \cite{Ham97}, that  $A$ is a subgroup if $|A|\ge 3.$ An  example  given in \cite{HamLlaLop} shows that this conclusion may fail in the non-abelian case.

\section{Arc-transitive digraphs}

Let $a\in V.$ The {\em twin class} of $a$ is $W_a=\{x: \Gamma (a)=\Gamma (x)\}.$ We need the following lemma.
%We also put $\lambda(a)=|\Lambda(a)|.$
%If this value is constant we shall denote it by $\lambda _{\Gamma}$ or simply by $\lambda,$
%when the context is clear.
\begin{lemma}
Let  $\Gamma=(V,E)$ be   vertex-transitive digraph and let  $v$ be an element of $V$.
Then  $ \{W_a:\ a\in V\}$ is a partition of $V$. Moreover $\Gamma^-(v)=\bigcup _{a\in \Gamma^-(v)} W_a.$
 Also
$|W_v|$
 divides
both $|V|$ and
$d(\Gamma).$ In particular, $\Gamma$ is irreducible if $\gcd(|V|,d(\Gamma))=1.$
\label{twinn}
\end{lemma}

\proof
 Assume that there is a $x$ with $x\in W_a \cap W_b,$ for some $a,b\in V.$
 We have $\Gamma (x)=\Gamma (a)=\Gamma (b).$ In particular, $a\in W_b$ and hence
 $W_a=W_b.$ It follows that $\{W_a:\ a\in V\}$ is a partition of $V$.
 In particular, $|W_v| $ divides $|V|.$

Assume that   $x\in \Gamma ^- (v)$ and let $y$ be an element of $W_x.$
We have $v\in \Gamma (x)$ and $\Gamma (y)=\Gamma (x).$ Thus $y\in \Gamma ^- (v).$
It follows that $\{W_a; a\in \Gamma^-(v)\}$ is a  partition of $\Gamma^-(v)$.
Take an arbitrary $x\in V$ and an automorphism $f$ such that $f(v)=x.$
Clearly $f(W_v)=W_x$. It follows that $|W_v|=|W_x|,$ for every $x\in V.$
Therefore,  $|W_v|$ divides both $|V|$ and
$d(\Gamma).$ \qed

%
%
%   EL LEMMA 4
%
%\begin{lemma}\label{connectedatom} Let $\Gamma =(V,E)$ be a $2$-separable finite  diconnectedstrong connected
  %arc-transitive digraph that has the intersection property of $2$-atoms and let $A$ be a $2$-atom of $\Gamma$ with $|A|\ge 3$. Then $\Gamma [A]$ is diconnected.
%\end{lemma}
%\proof
The next result gives a good description for vosperian arc-transitive digraphs.

\begin{theorem} A finite   strongly connected
  arc-transitive  digraph  $\Gamma =(V,E)$  of  degree $d\le |V|-4$ with $d\notin
\{1,2,4,6\}$
is vosperian if and only if it is irreducible.

\label{arctransitive} \end{theorem}

\proof
As we observed in the Introduction a reducible digraph is non-vosperian.
 Assume that $\Gamma$ is  an irreducible arc-transitive digraph which is non-vosperian. We prove the following points:

(i)  No  block  is a $2$-fragment.

Suppose to the contrary that some block $B$ is a $2$-fragment.
Take $a\in B$ and  an arc $(a,b)$ with $b\notin B.$ The set $B$ must
be an independent set,
otherwise there would exist an arc $(c,d)$ inside $B$. Since $\Gamma$ is arc-transitive, there is an automorphism $f$ with $f(c)=a,$
and $f(d)=b.$ Thus $B\cap f(B)\neq \emptyset$ but $B\neq f(B)$, contradicting the definition of  a block. Hence, $\Gamma (x)\subset
\partial (B),$ for every $x\in B.$ Since $|\Gamma (x)|\le |\partial (B)|\le\kappa _2(\Gamma)\le d,$ we have
$\Gamma (x)= \partial (B),$ for every $x\in B.$ In particular $\Gamma$
has a twin pair, a contradiction.

(ii)  No  block  is  a negative $2$-fragment.
The proof is similar to the proof of (i).

(iii) Every vertex of $V$ is contained in at least three
distinct $2$-atoms or there are distinct $2$-atoms $M,N$ with $|M\cap N|\ge 2.$

Suppose on the contrary, that two distinct $2$--atoms intersect in
at most one point, so that every pair of vertices is contained in at
most one $2$--atom, and  that every vertex is in at most two
distinct $2$--atoms.

Take a $2$-atom $A.$ Suppose that $|A|=2$. Since $|\partial
A|=\kappa_2(A)\le d$ (otherwise $\Gamma$ is vosperian) and the
digraph  has no twin pairs, $\Gamma [A]$ contains an arc $(x,y)$.
Since the degree $d\ge 3$, the automorphisms sending $(x,y)$ to
$(x,y')$ for each $y'\in \Gamma (x)$ provide at least three distinct
$2$--atoms containing $x$. Hence we may assume   that $|A|> 2$.

Let us show that $\Gamma [A]$ is a strongly connected vertex transitive
digraph. By Lemma \ref{katom},  $\Gamma [A]$ contains a circuit, say
$C=[a_1,a_2, \cdots, a_j].$ Take an arbitrary
 vertex $v\in A.$  By Lemma \ref{katom}, there is $w\in A\cap \Gamma^-(v).$
 Take an automorphism $f$ with $f(a_1,a_2)=(w,v)$. Since every pair
 of vertices is contained in a unique atom, we have $f(A)=A$.
 Therefore
$f(C)$ is a circuit in $\Gamma [A]$ containing $v.$ Thus every arc
of $\Gamma [A]$ is contained in a circuit of length at least two.
If $\Gamma [A]$ is not strongly connected, it contains a strongly connected component
$K$ which is a sink, which must have cardinality at least two. We
have $\partial (K)\subset \partial (A).$
 Hence, by the minimality of a $2$-atom  $K=A.$
 It follows that $\Gamma [A]$ is  strongly connected.
Moreover, the same argument shows that the arbitrary vertex
 $v$ can be sent to $a_2$ by an automorphism which leaves $\Gamma
 [A]$ invariant, so that this induced subgraph is vertex transitive and hence
  a regular digraph.
We shall denote the degree of $\Gamma
[A]$ by $r.$

There is an element $c\in A$ with an element $c'\in \Gamma (c)\cap \partial (A)$. By Lemma \ref{katom}, there is an arc
  $(b,c)$ contained in $A.$
  For every arc
$(x,y)\in E,$ there is an automorphism $f$
with $f (b)=x$ and $f (c)=y.$ In particular $f (A)$ is a $2$-atom containing
$\{x,y\}.$
By our hypothesis, such a $2$-atom  is  unique.
We shall denote this atom by $A_{xy}.$
Since $c$ is contained in at most two distinct $2$-atoms, $\Gamma (c)=(\Gamma(c)\cap A)\cup
  (\Gamma(c)\cap A_{cc'}).$
In particular,
$d=2r.$
Let $a$ be an element in $A\setminus\{b,c\}$ and let $a',b'$ be vertices in the boundary of $A$ such that $a'\in \Gamma (a)$ and $b'\in \Gamma (b)$.
  We have $\partial (A)\supset (\Gamma (a)\cap A_{aa'})\cup (\Gamma (b)\cap
A_{bb'})\cup (\Gamma (c)\cap A_{cc'}).$
  Thus we have using our hypothesis on atoms intersection,
   \begin{align*}
  |\partial (A)|&\ge |\Gamma
(a)\cap A_{aa'}|+|(\Gamma (b)\cap A_{bb'})\setminus A_{aa'}|+|(\Gamma (c)\cap
A_{cc'})\setminus (A_{aa'}\cup A_{bb'})|\\
&\ge |\Gamma
(a)\cap A_{aa'}|+|(\Gamma (b)\cap A_{bb'}|-1+|\Gamma (c)\cap
A_{cc'}|-2\\
&= r+r-1+r-2=
  3d/2-3>d,
   \end{align*}
   since $d>6,$ this implies a contradiction.

(iv) Every vertex of $V$ is contained in at least three
distinct negative $2$-atoms or there are distinct negative $2$-atoms $M,N$ with $|M\cap N|\ge 2.$
   The proof is exactly the same as for (iii).

By Theorem \ref{asterix}, $\Gamma$ has a $2$-atom with size $2$ or a negative
$2$-atom with size $2.$
The two cases are similar, we choose the first one. Since $\Gamma$ has
no twins, the $2$-atom
has the form $A=\{a,b\},$ with $(a,b)\in E.$

Let $v$ be the only element of
$\Gamma (b)\setminus \Gamma (a)$ and take $u\in \Gamma (b)\cap \Gamma
(a)$.
Now $\{b,v\}$ is a $2$-atom. It follows that $u\in \Gamma (v)$. Also
$\{b,u\}$ is a $2$-atom. It follows that $v\in \Gamma (u)$, a contradiction in the anti-symmetric case. Assume now that   $\Gamma$ is symmetric. A
similar reasoning proves that, $u\in\Gamma
(z)$ where $z$ is the only element of $\Gamma (a)\setminus
\Gamma (b)$. Since $\{a,u\}$ is a  $2$-atoms, it follows that
$|\Gamma(u)\cap(\Gamma (b)\cap \Gamma(a))|\ge d-4$. Thus, $|\Gamma(u)\cap
\Gamma (\{a,b\})|=d$ and hence, $\Gamma (u)\cap \nabla
(\{a,b\})=\emptyset$, a contradiction.\qed

\begin{corollary}
Let $\Gamma$ be   a  strongly connected arc-transitive  digraph  of degree $d\notin
\{2,4,6\}$  such that $\gcd(|V|,d)=1.$
Then  $\Gamma$ is vosperian.
\label{superung1}
\end{corollary}
\proof
By Lemma
\ref{twinn},
$\Gamma$ is irreducible.
 The result follows now by Theorem \ref{arctransitive}.\qed

\section{Vosperian  Cayley  digraphs}

In this section, we investigate  vosperian   Cayley digraphs.

 Let $X$ be a subset of a group $G.$ We shall  write $\tilde{X}=X\cup \{1\}.$
 We shall say that $X$ is  a {\em right $r$-coprogression},
if $G\setminus X=\{a, ra, \ldots , r^ja\},$ for some $a\in G$.
Similarly we define a {\em left $r$-coprogression}.
Notice that $S$ is  a {right} $r$-coprogression if and only if
$S$ is a left $a^{-1}ra$-coprogression.

In particular,  if $\tilde{S}$ is a right $r$-coprogression and $|S|\le |G|-4$ then \begin{equation}\label{coprogression}|\Gamma(\{1,r\})\cup \{1,r\}|=|\{1,r,\}\tilde{S}|=|S|+2,\end{equation} and
$\kappa_2(\Gamma )\le |S|.$

  \begin{theorem}
Let $S$ be a   generating subset of a finite group $G$, with $1\notin S.$
 Then $\Gamma =\mbox{Cay} (G,S)$ is  non-vosperian if and only if one of the following holds:

  \begin{enumerate}

    \item [(i)] There are a subgroup $H$  of $G$ with $|H|\ge 2,$ and  an element $a\in G$ such that
$|(H\cup Ha)\tilde{S}|\le \min (|G|-2,|H\cup Ha|+|{S}|).$
    \item  [(ii)] There are a subgroup $H$  of $G$ with $|H|\ge 2,$ and  an element $a\in G$ such that
$ |\tilde{S}(H\cup aH)|\le \min (|G|-2,|H\cup Ha|+|{S}|).$
    \item  [(iii)]  There is a $r\in G\setminus \{1\}$ such that $\tilde{S}$ is  a right $r$-coprogression and $|G|-4\ge |S|.$
  \end{enumerate}

\label{vosperg}
\end{theorem}

\proof
Let us first prove the sufficiency.
%Notice that  $\Gamma$ is non-vosperian if and only if $\Gamma^{-1}$ is non-vosperian.
If there exist a subgroup $H$  with $|H\ge 2$ and an element $a\in G$ such that  (i) or (ii)  holds,  then $\Gamma $
 is $2$-separable and $\kappa_2(\Gamma )\le |S|$.
Similarly, if there exists a $r\in G\setminus \{1\}$, such that $\tilde{S}$ is a right $r$-coprogression, then
$|\Gamma(\{1,r\})\cup \{1,r\}|=|\{1,r,\}\tilde{S}|=|S|+2,$  and hence
 $\Gamma $ is $2$-separable and $\kappa_2(\Gamma )\le |S|$.   So each of these conditions is a necessary one.

%Let us now prove the sufficiency:

Suppose now that $\Gamma$ is non-vosperian.
%In particular $\Gamma$ is $2$-separable and hence and $|G|-4\ge |S|$.
%, thus  by Lemma \ref{equal_kappa} and by Lemma \ref{cvosper}, ${\Gamma} $ and $\Gamma^-$ are  $2$-separable and $\kappa _2({\Gamma} )=\kappa
%_2({\Gamma}^- )\le |S|.$
% Let $ M$ be a $2$-atom of $ \Gamma$  and
%let $N$ be a negative $2$-atom of $ \Gamma$  such that  $1\in M\capN.$
% It would be enough to consider the case where $|M|\le |N|,$ the other cases reduces to this one by replacing $\Gamma$ by $\Gamma ^-.$
  By Corollary \ref{vospercorog}, there are a subgroup $H$ and an element $a\in G$
such that $H\cup Ha$ is a $2$-fragment or a negative $2$-fragment. The two cases are similar (and equivalent up to duality), so we shall consider
 only the case where $H\cup Ha$ is a $2$-fragment.
Assume first that  $|H|\ge 2.$
Since $H\cup Ha$ is a $2$-fragment, we have $|(H\cup Ha)\tilde{S}|=|\Gamma (H\cup Ha)|\le |G|-2.$ We have also
$|H\cup Ha|+|S|\ge |H\cup Ha|+\kappa_2(\Gamma)=|(H\cup Ha)\tilde{S}|.$  Therefore   (i) holds.

Assume now that $|H|=1.$
 Let $M=\{1,r\}$. By the definition of a $2$-atom we have $|\{1,r\}\tilde{S}|=\kappa_2(\Gamma)+2\le |{S}|+2.$
 We may assume that $\kappa_2(\Gamma)=|S|.$ Otherwise $r\tilde{S}=\tilde{S},$ and hence $\tilde{S}$ is a union of right $\subgp{r}$-cosets, and thus Condition (i) holds with $H=\subgp{r}$.
  Hence,

\begin{equation}\label{eq2g}
\kappa_2(\Gamma)=|S| \ \text{and thus } \ |\{1,r\}\tilde{S}|=|{S}|+2.
\end{equation}

 Let $K$ be the cyclic subgroup generated by $r$.
If $G=K,$ then by (\ref{eq2g}), $\tilde{S}$ is a progression and also a coprogression since $r$ generates $G.$
So we may assume that $K\neq G.$

Let
$\tilde{S}=S_1\cup \cdots\cup S_j$ be the partition of $\tilde{S}$ induced by the partition of $G$ into right $K$-cosets. We have $j\ge 2$ since $K\neq G$. We shall assume that $|S_1|\le
\cdots \le |S_j|$. We must have $|S_2|=|K|,$ since otherwise $|\{1,r\}{S_i}|\ge |S_i|+1,$ for all
$1\le i \le 2$. It would follows that $|\{1,r\}\tilde{S}|\ge |S_1|+1+ |S_2|+1 +\sum_{i\ge 3}|S_i|\ge |\tilde{S}|+2,$
contradicting (\ref{eq2g}). It follows also that $|\{1,r\}{S_1}|= |S_1|+1.$ Since $r$ generates $K$, Lemma \ref{AP} implies that
$S_1$ is a right progression with ratio $r$. Note that $KS_1\setminus S_1$ is also a right progression.

{\bf Subcase 2.1} $j|K|=|G|.$ In this case $G\setminus \tilde{S}=KS_1\setminus S_1$ is also a right progression.

{\bf Subcase 2.2} $j|K|\le |G|-|K|.$ We  have
$(j-1)|K|+|S_1|-1=|\tilde{S}|-1=\kappa _2(\Gamma)\le |K\tilde{S}|-|K|=(j-1)|K|.$ Thus $|S_1|=1.$
In particular, $K$ is clearly a $2$-fragment and (i) holds with $H=K$ and $a=1$. \qed

 %{\bf Case} 4: $|N|=2.$  Since this case reduces to Case 3, by replacing $S$ with $S^{-1},$ the proof is completed.

  %%%%%%%%%%%%%%%%%%%

\section{Superconnected Cayley digraphs}

In this section we characterize irreducible  superconnected Cayley digraphs.
As we have seen, these Cayley digraphs are defined by an aperiodic subsets.

\begin{theorem}
An irreducible  strongly connected Cayley digraph $\Gamma=Cay(G,S)$  on a finite group $G$
is superconnected
if and only if one of the following conditions holds:
\begin{itemize}
 \item[(i)]  $\Gamma$ is vosperian,
  \item[(ii)] for some $r\in G,$  $\tilde{S}$ is a right $r$-coprogression with $ r^{-1}\notin S$ and $|G|-4\ge |S|$.
 \end{itemize}
  %Moreover $a\in H$, if $G$ is abelian.
\label{superg}
\end{theorem}

\proof
Take  a  subset $S$ of  $G$
such  that  $\Gamma=\mbox{Cay} (G,S).$
Since $\Gamma$ is  strongly connected, $S$ must be a generating subset of $G.$

 Condition (i) implies obviously that $\Gamma$ is superconnected.  Assume that (ii) holds
and let $P=\subgp{r}\cap \tilde{S}.$ By Lemma \ref{twincayley},  $S$ is
aperiodic, and thus $|P|\ge 2.$ Note that, since $\tilde{S}$ is a coprogression $\kappa (\Gamma)=|S|$. Take now any fragment $F$ of $\Gamma$. By a left-translation we can assume that
$1\in F.$ We must have $F \subset H,$ since otherwise $|FH|\ge
2|H|,$ and hence by Lemma \ref{prehistorical}, $$|F\tilde{S}|\ge
|F(\tilde{S}\setminus H)|=|(F)H(\tilde{S}\setminus
H)|=|(FH)(\tilde{S}\setminus H)|=|G|,$$ a contradiction. Now we have
$|F|+|\tilde{S}|-1=|F\tilde{S}|=|G\setminus
H|+|FP|=|\tilde{S}|-|P|+|FP|.$ In particular, $|H|-1\ge
|FP|=|F|+|P|-1.$ By Lemma \ref{AP}, $F$ is a $r$-progression. By a
left-translation, we may assume that $F=\{1,r, \ldots ,
r^{|F|-1}\}.$ Since $|H|-1\ge |FP|,$ we have $\Gamma  ( r^{|F|-1})=
\partial (F).$ Thus, $\mbox{Cay}(G,S)$ is  superconnected.

Let us now prove the necessity.

Suppose that  (i) and (ii) do not hold. Since $\Gamma$ is non-vosperian, by Theorem \ref{vosperg} we are in one of the following cases:

{\bf Case} 1. There are a subgroup $H$ and an element $a\in G$ with $|H|\ge 2$ such that
   $|G|-2\ge|(H\cup Ha) \tilde{S}|=|H\cup Ha|+| \tilde{S}|-1$.

   Put $T=(H\cup Ha) \tilde{S}\setminus (H\cup Ha).$ Clearly $T$ is a minimum cutset  with $HT=T.$
   We can not have $T=\Gamma (x),$ for some $x,$ otherwise $HxS=xS,$ and $S$ would be right-periodic, a contradiction.
   Also, we can not have $T=\Gamma^- (y),$ for some $y,$ otherwise $HyS^{-1}=yS^{-1},$
   and $S^{-1}$ would be right-periodic, a contradiction.
   Thus $\Gamma$  is non superconnected.

{\bf Case} 2.  There are a subgroup $H$ and an element $a\in G$ with $|H|\ge 2$ such that
    $|G|-2\ge|(H\cup Ha) \tilde{S}^{-1}|=|H\cup Ha|+| \tilde{S}|-1$. This case is similar to the previous one.

 {\bf Case} 3. For some $r\in G$,  $\tilde{S}$ is a right $r$-coprogression and $|S|\le |G|-4$.  Put $H=\subgp{r}$ and $G\setminus \tilde{S}\subset Ha$ for some $a\in G$.

{\bf Subcase} 3.1 $a\notin H$.
Since $\tilde{S},$ is a coprogression, we have $H\subset \tilde{S}.$
Since  $(a^{-1}Ha)\cap Ha=\emptyset,$ we have  $a^{-1}Ha\subset \tilde{S}.$
Take an arbitrary $h\in H\setminus \{1\}.$ We have $$\Gamma (ha)=haS\supset ha(a^{-1}Ha\setminus \{1\})=Ha\setminus \{ha\}.$$
It follows that $\Gamma[Ha]$ is a complete symmetric digraph.
Since $H \subset \tilde{S},$ we have $r, r^{-1} \in S$. Thus, $\Gamma[\{1,r\}]$ is a complete symmetric digraph.
Hence $V\setminus \partial (\{1,r\})$ has exactly two  strongly connected components each of them has size $\ge 2.$

Let $T=\partial (\{1,r\}).$ Since $|\{1,r\}\tilde{S}|=|S|+2,$ we have $|T|=|S|.$
By the observation made above $T\neq \Gamma (x),$ for every $x\in G.$ Therefore,  $\Gamma$ is non superconnected.

{\bf Subcase} 3.2 $a\in H$.  Since (ii) does not hold, we have $r, r^{-1} \in S$. Put
$P=\subgp{r}\cap \tilde{S}$.
 Let $F$ be any fragment with $1\in F.$ We must have $F \subset H,$ since otherwise by Lemma \ref{prehistorical}, $$|F\tilde{S}|\ge
|F(\tilde{S}\setminus H)|=|(F)H(\tilde{S}\setminus H)|=|(FH)(\tilde{S}\setminus H)|=|G|,$$  a contradiction.

Now we have $|F|+|\tilde{S}|-1=|F\tilde{S}|=|G\setminus H|+|FP|=|\tilde{S}|-|P|+|FP|.$
In particular, $|H|-1\ge |FP|=|F|+|P|-1.$ By Lemma \ref{AP}, $F$ is a $r$-progression.
Since $\{r,-r\}\subset S,$ the set $F$ is a  strongly connected subset. Similarly $G\setminus \Gamma (F)$
is a  strongly connected  subset. Therefore $\Gamma$ is non superconnected.\qed

\begin{corollary}
Let $S$ be   an aperiodic  generating subset of a finite group $G$,
with $1\notin S,$   and  $ |S|\le |G|/2$. Then  $\mbox{Cay} (G,S)$
is   superconnected if and only if  one of the following holds:

\begin{itemize}
  \item  $\mbox{Cay} (G,S)$ is vosperian,
  \item  There is a $r\in G$ such that $G$ is a cyclic group and $S=\{r,r^2,\ldots ,r^{|S|}\}.$
\end{itemize}
 \label{superung2}
\end{corollary}

\proof
Let us see that if $ |S|\le |G|/2$ holds then a coprogression is also a progression. If $G\setminus \tilde{S}\subset \subgp{r}$
then $G=\subgp{r}$. Since otherwise, $|\tilde{S}|\ge |G|-|\subgp{a}|+2\ge \frac{|G|}2+2,$
a contradiction. The result follows now by Theorem \ref{superg}.\qed

\noindent{\bf Acknowledgment}

Research supported by the Ministry of Science and Innovation, Spain under project MTM2008-06620-C03-01/MTM. Research done when the last author was \mbox{visiting} Universit\'e Pierre et Marie Curie, E. Combinatoire, Paris, supported by the Ministry of Science and Innovation, Spain  under the National Mobility Programme of Human Resources, Spanish National Programme I-D-I 2008--2011.

\bibliographystyle{model1a-num-names}
%\bibsection

\end{document}